\documentclass[12pt,reqno]{amsart}

\usepackage{amssymb,latexsym}

\usepackage{enumerate}
\allowdisplaybreaks
\usepackage[french,english]{babel}
\usepackage{amsmath}
\usepackage{graphicx}
\usepackage{amssymb}
\usepackage{bbm}
\usepackage{amsthm,mathtools}
\usepackage{ulem}
\usepackage{geometry}
\usepackage{tikz-cd}
\usepackage{mathrsfs}
\usepackage[colorinlistoftodos]{todonotes}
\usepackage{enumitem}
\usepackage{verbatim}
\usepackage[foot]{amsaddr}
\usepackage{dsfont}
\usepackage[T1]{fontenc} 
\usepackage{cite}

\makeatletter

\@namedef{subjclassname@2010}{
	
	\textup{2020} Mathematics Subject Classification}

\makeatother
\newtheorem{thm}{Theorem}[section]
\newtheorem*{thm*}{Theorem}

\newtheorem{prop}{Proposition}[section]
\newtheorem{lem}[thm]{Lemma}

\theoremstyle{definition}

\newtheorem*{xrem}{Remark}
\numberwithin{equation}{section}

\newcommand{\inv}{^{-1}}

\newcommand{\mbd}{\mathbb{D}}
\newcommand{\mbk}{\mathbb{K}}

\newcommand{\mbz}{\mathbb{Z}}
\newcommand{\mbr}{\mathbb{R}}

\newcommand{\mbf}{\mathbb{F}}

\newcommand{\mbn}{\mathbb{N}}

\newcommand{\mbq}{\mathbb{Q}}

\newcommand{\mch}{\mathcal{H}}
\newcommand{\mca}{\mathcal{A}}

\newcommand{\mce}{\mathcal{E}}
\newcommand{\mcf}{\mathcal{F}}
\newcommand{\mcg}{\mathcal{G}}
\newcommand{\mcj}{\mathcal{J}}

\newcommand{\mcn}{\mathcal{N}}
\newcommand{\mcp}{\mathcal{P}}

\newcommand{\mcm}{\mathcal{M}}

\newcommand{\mcz}{\mathcal{Z}}

\usepackage{hyperref}
\hypersetup{colorlinks=true,linkcolor=blue,anchorcolor=blue,citecolor=blue}      

\usepackage{color}
\newcommand{\newabstract}[1]{%
	\par\bigskip
	\csname otherlanguage*\endcsname{#1}%
	\csname captions#1\endcsname
	\item[\hskip\labelsep\scshape\abstractname.]
}

\begin{document}

	\baselineskip=17pt

	\title[Extreme values of derivatives of the Dedekind zeta function of a cyclotomic field]{Extreme values of derivatives of the Dedekind zeta function of a cyclotomic field}

    \author{Zhonghua Li\textsuperscript{1}}
    \author{Yutong Song\textsuperscript{1,2}}
    \author{Qiyu Yang\textsuperscript{3}}
    \author{Shengbo Zhao\textsuperscript{1}}
    \address{1.School of Mathematical Sciences, Key Laboratory of Intelligent Computing and Applications(Ministry of Education), Tongji University, Shanghai 200092, P. R. China}
    \address{2.Graz University of Technology, Institute of Analysis and Number Theory, Steyrergasse 30/II, 8010 Graz, Austria}
    \address{3.School of Mathematics and Statistics, Henan Normal University, Xinxiang 453007, CHINA}
    \email{zhonghua\_li@tongji.edu.cn}
    \email{99yutongsong@gmail.com} 
    \email{qyyang.must@gmail.com}
	\email{shengbozhao@hotmail.com}

	\begin{abstract} 
	    In this paper, we establish a lower bound for maximum of derivatives of the Dedekind zeta function of a cyclotomic field on the critical line. Employing a double-version convolution formula and combined with special GCD sums, our result generalizes the work of Bondarenko et al. (2023). We also establish a lower bound via the resonance method when the real part is near the critical line, both of the above results refine parts of Yang's work (2022).  
	\end{abstract}
	
    \keywords{The Dedekind zeta function, extreme values, GCD sums, resonance method. }
	
	\subjclass[2020]{Primary 11M06, 11M26, 11N37.}
	
	\maketitle
\section{Introduction}
\label{Introduction}
Extreme values have long been central objects of study in analytic number theory, as they play a crucial role in understanding the Riemann zeta function $\zeta(s)$ and Dirichlet L-functions $L(s,\chi)$. Such quantities are instrumental in investigating the distribution of values of these functions. In \cite{soundararajan2008extreme}, Soundararajan introduced a powerful technique, known as the resonance method, for studying extreme values for various families of $L(s, \chi)$. When applied to $\zeta(s)$, he showed that for sufficiently large $T$,
$$
\max_{T \le t \le 2T} \Big|\zeta\Big(\frac{1}{2}+it\Big)\Big| \ge \exp\bigg((1 + o(1))\sqrt{\frac{\log T}{\log_2 T}}\bigg).
$$
Here and throughout, we write $\log_j$ for the $j$-th iterated logarithm.
\par
Combining ideas from Aistleitner \cite{aistleitner2016lower}, Bondarenko and Seip \cite{bondarenko2017large,bondarenko2018argument} refined Soundararajan's method and proved that for sufficiently large $T$,
$$
\max_{\sqrt{T} \le t \le T} \Big| \zeta\Big( \frac{1}{2} + it \Big) \Big| \geq \exp\bigg( (1+o(1)) \sqrt{ \frac{\log T \log_3 T}{\log_2 T} } \bigg).
$$
Later, the leading term $(1+o(1))$ was improved to $(\sqrt{2}+o(1))$ by de la Bretèche and Tenenbaum \cite{tenen2019galsum}, by invoking optimized GCD sums.
\par
Recently, Yang \cite{yang2022extreme,yang2024extreme} extended these ideas to derivatives of $\zeta(s)$, obtaining analogous lower bounds in the critical strip. Denote by $\mbn$ the set of all non-negative integers. When $T \to \infty$, he established the following results concerning GCD sums for $\ell \in \mbn$:
$$\max _{t \in[0, T]}\Big|\zeta^{(\ell)}\Big(\frac{1}{2}+i t\Big)\Big| \ge\exp \bigg(\sqrt{2} \cdot\sqrt{\frac{\log T \log_3 T}{\log_2T}}\bigg)$$
and
$$\max _{t \in[0, T]}\Big|\zeta^{(\ell)}(\sigma+i t)\Big| \ge \exp \bigg(\frac{c}{1-\sigma} \cdot \frac{(\log T)^{1-\sigma}}{(\log \log T)^\sigma}\bigg),$$
where $\sigma \in (1/2, 1)$ and $c>0$ is an absolute constant.
\par
Here and throughout this paper, let $\mbk = \mbq(\omega_d)$, where $\omega_d$ is a $d$th root of unity for $d \ge 3$. One of the main motivations for studying the Dedekind zeta function $\zeta_\mbk(s)$ is that it satisfies the following classical factorization
\begin{equation*}
    \zeta_\mbk(s)=\zeta(s)\prod_{\chi\neq\chi_{0}(\text{mod }d)}L(s,\chi^{\ast}),
\end{equation*}
where the product runs over all non-principal characters $\chi$ modulo $d$, and $\chi^{\ast}$ is the character that induces $\chi$ if $\chi$ is not primitive and $\chi^{\ast} = \chi$ otherwise. Thus, sharper extreme values of $L(s, \chi)$ may be obtained by studying $\zeta_\mbk(s)$. For a general field $\mbf$, Li \cite{li2018extreme} showed that there exist arbitrarily large $t$ such that
$$\Big|\zeta_\mbf\Big(\frac{1}{2}+i t\Big)\Big| \ge \exp \bigg(c \sqrt{\frac{\log t \log_3 t}{\log_2 t}}\bigg),$$
the magnitude of $c$ depends on whether $\mbf/\mbq$ is a Galois extension. 
\par
Subsequently, Bondarenko et al.  \cite{bondarenko2023dichotomy} optimized Li's work and showed that for sufficiently large $T$,
$$\max _{t \in[0, T]}\Big|\zeta_\mbk\Big(\frac{1}{2}+i t\Big)\Big| \ge \exp \bigg(\Big(\sqrt{\phi(d)} +o(1)\Big) \sqrt{\frac{\log T \log_3 T}{\log_2T}}\bigg)$$
uniformly for $d \ll (\log_2T)^A$ with an arbitrary positive number $A.$ Here, $\phi(n)$ is the Euler's totient function. This conclusion reveals that it is possible to obtain sharper $\Omega$-results for
$\zeta(s)$ and $L(s,\chi)$. Moreover, the best known bound for $L(s, \chi)$
$$
\max _{T \le t \le 2T}\Big| L\Big(\frac{1}{2} + i t, \chi \Big)\Big| \ge \exp \bigg((1+o(1)) \sqrt{\frac{\log T}{\log_2 T}}\bigg)
$$
was established by Soundararajan in \cite{soundararajan2008extreme}. 
\par
Suppose that $\mcm \subset \mbn$ is a finite set, and let $\sigma \in (0, 1].$ Define the GCD sums $S_{\sigma}(\mcm)$ as 
$$
S_{\sigma}(\mcm) \coloneqq \sum_{m,n\in\mcm} \frac{(m,n)^{\sigma}}{[m,n]^{\sigma}},
$$
where $[m,n]$ and $(m,n)$ denote the least common multiple and the greatest common divisor of $m$ and $n$ respectively. G{\'a}l \cite{gal1949theorem} derived that $\sup_{|\mcm| = N} S_1(\mcm)/|\mcm| \asymp (\log_2 N)^2$ for sufficiently large $N$.
Later, de la Bretèche and Tenenbaum \cite{tenen2019galsum} investigated $S_{1/2}(\mcm)$ and obtained 
$$\sup_{|\mcm| = N} \frac{S_{1/2}(\mcm)}{|\mcm|} = \exp\bigg( \Big(2\sqrt{2} + o(1)\Big) \sqrt{\frac{\log N \log_3 N}{\log_2 N}} \bigg)$$
as $N \to \infty$. For $1/2 < \sigma < 1,$ Aistleitner, Berkes, and Seip \cite{aistleitner2015gcd} proved that
$$
\exp\bigg(c_{\sigma}\cdot\frac{(\log N)^{1 - \sigma}}{(\log_2 N)^{\sigma}}\bigg) \ll \sup_{|\mcm| = N}\frac{S_{\sigma}(\mcm)}{|\mcm|} \ll \exp\bigg(C_{\sigma}\cdot\frac{(\log N)^{1 - \sigma}}{(\log_2 N)^{\sigma}}\bigg),
$$
where $c_{\sigma}$ and $C_{\sigma}$ are positive constants only depending on $\sigma$. They obtained the above bounds by utilizing the trigonometric polynomials on infinite-dimensional polydisc $\mbd^{\infty}$. Notably, their result is almost optimal. GCD sums are powerful tools for studying extreme values of functions related to $\zeta(s)$. Several studies employing GCD sums, including those by Bondarenko and Seip \cite{bondarenko2017large}, de la Bretèche and Tenenbaum \cite{tenen2019galsum}, and Yang \cite{yang2022extreme}, have yielded significant results. For more results and details about GCD sums, we recommend \cite{bondarenko2016gal,bondarenko2015gcd,bondarenko2018note,dong2021large,lewko2017refinements,yang2022note} and the references therein.
\par
Motivated by 
 \cite{bondarenko2023dichotomy,yang2022extreme}, we study extreme values of derivatives of $\zeta_\mbk(s)$ in the critical strip. We know that
$$\zeta_{\mbk}(s) = \sum_{n=1}^{\infty}\frac{a_{\mbk}(n)}{n^s},$$
where $\Re(s) > 1$ and $a_{\mbk}(n)$ is a multiplicative function. Let $p$ be a prime number. By \cite{washington1997introduction}, we have
\begin{equation*}
    a_{\mbk}(p^k) = 
    \begin{cases}
        N(\frac{k}{f},r), &\, \operatorname{if} f \mid k, \\
        0, &\, \operatorname{otherwise},
     \end{cases}
\end{equation*}
where $f$ denotes the multiplicative order of $p$ in $(\mbz/d\mbz)^*$, and $N(k/f,r)$ is the number of ways to write $k/f$ as a sum of $r$ non-negative integers.
\par
Henceforth, we assume that $\ell \in \mbn$. Our main results, Theorems \ref{thm1.1} and \ref{thm1.2}, establish lower bounds for the maximum of derivatives of $\zeta_\mbk(s)$ both on and near the critical line. Specifically, Theorem \ref{thm1.1} shows that the lower bound on the critical line attains the same quality as that for $\zeta_\mbk(s)$ itself, while Theorem \ref{thm1.2} extends the estimate to a suitable range.
\begin{thm}
    \label{thm1.1}
   Let $A$ be an arbitrary positive number. If $T$ is sufficiently large, then uniformly for $d \ll (\log_2 T)^A,$ we have
     $$
\max_{t\in[0,T]}\Big|\zeta_{\mbk}^{(\ell)}\Big(\frac{1}{2}+it\Big)\Big| \ge \exp\bigg(\Big(\sqrt{\phi(d)} + o(1)\Big)\sqrt{\frac{\log T\log_3 T}{\log_2 T}}\bigg).
     $$
\end{thm}
\par
Bondarenko et al. stated in \cite{bondarenko2023dichotomy} that it is likely that the methods of de la Bretèche and Tenenbaum \cite{tenen2019galsum} can improve the exponent on
the right-hand side of the above by a factor of $\sqrt{2}$ when $\ell = 0$. However, since $a_\mbk(n)$ is not a completely multiplicative function, we have to employ a certain function $a_\mbk^\prime(n)$ to establish the lower bound for GCD sums. This procedure results in the loss of a factor of $\sqrt{2}$. See the Appendix for the details.
\par
Theorem \ref{thm1.1} refines certain results of \cite{bondarenko2023dichotomy,yang2022extreme}. In comparison with \cite{bondarenko2023dichotomy}, we derive a lower bound for maximum of derivatives of $\zeta_\mbk(s)$ on the critical line. When $\ell = 0,$ it reduces to the case of $\beta = 0$ in \cite[Theorem 1]{bondarenko2023dichotomy}. This indicates that the lower bound for the maximum of derivatives of $\zeta_\mbk(s)$ on the critical line can reach the same order as that of the original Dedekind zeta function $\zeta_\mbk(s)$. Comparing with \cite{bondarenko2023dichotomy}, in addition to the above refinement, we also enhance the factor $\sqrt{\phi(d)}$ in the exponent on the right-hand side of \cite[Eq. (2)]{bondarenko2023dichotomy}. Camparing with \cite{yang2022extreme}, we obtain a result that generalizes $\zeta(s)$. This result corresponds to the case of $\beta = 0$ in \cite[Theorem 2.(A)]{yang2022extreme}. It should be noted that when $d = 1,2,$ then $\mbk=\mbq$ is the rational number field with $\phi(d)=1$. In this case, Theorem \ref{thm1.1} is weaker than \cite[Theorem 2]{yang2022extreme}. Therefore, we focus on the cases where $d\ge 3$, that is, when $\mbk$ is not the trivial rational number field $\mbq$. 
\par
When the real part $\sigma$ is near the critical line, we state Theorem \ref{thm1.2} below.
\begin{thm}
    \label{thm1.2}
    Let $A$ be an arbitrary positive number. Let $\sigma > 0$ be a real number and $T>0$ sufficiently large in the range 
    $$\frac{1}{2} \le \sigma \le \frac{1}{2} + \frac{1}{\log_2 T}.$$
    Then uniformly for $d \ll (\log_2 T)^A,$ we have
    $$
    \max_{t\in[0,T]}\Big|\zeta_{\mbk}^{(\ell)}(\sigma+it)\Big|\ge \exp\bigg(\Big(\sqrt{\frac{\phi(d)}{e-1}}+o(1)\Big)\frac{(\log T)^{1-\sigma}(\log_3 T)^{\sigma}}{(\log_2 T)^{\sigma}}\bigg).
    $$
\end{thm}
\par
Theorem \ref{thm1.2} provides a lower bound for the maximum of derivatives of $\zeta_\mbk(s)$ near the critical line. Compared with \cite[Theorem 1]{bondarenko2023dichotomy}, our leading term in the exponent becomes smaller because the series $\sum p^{-2\sigma}$ can only obtain an upper bound similar to \cite[Lemma 9]{chirre2019extreme}. We can also apply the same method used in the proof of Theorem \ref{thm1.2} to study extreme values of derivatives of $\zeta_\mbk(s)$ on the critical line, but comparing with Theorem \ref{thm1.1}, the leading term in the exponent here is much smaller.
\par
To obtain Theorems \ref{thm1.1} and \ref{thm1.2}, we rely on the resonance method. This approach was originally proposed by Voronin \cite{voronin1988lower}, although his work did not attract much attention at the time. In 2008, Soundararajan \cite{soundararajan2008extreme} developed the resonance method to study extreme values of $\zeta(s)$ and $L(s, \chi)$. Subsequently, Aistleitner \cite{aistleitner2016lower} further refined Soundararajan's resonance method and improved the results of \cite{voronin1988lower} and \cite{hilberdink2009arithmetical}. The resonance method we employ is similar to that refined by Bondarenko and Seip \cite{bondarenko2017large,bondarenko2018argument}, which is following ideas of \cite{aistleitner2016lower}, and the slightly modified approach by Chirre \cite{chirre2019extreme}. For more results and details about the resonance method, we recommend \cite{aistleitner2019extreme,dong2023Onde,xumax2024extreme,qiyu2024large} and the references therein.
\par
Then, we introduce some notations. Let $A$ be an arbitrarily large positive number and $\varepsilon$ be small. We note that each occurrences of $A$ and $\varepsilon$ may represent different values. Let $\mbr$ denote the set of all real numbers. Finally, we denote the Fourier transform of a function $f \in L^1(\mbr)$ as
$$\widehat{f}(\xi)\coloneqq \int_\mbr f(x) e^{-i x \xi} \mathrm{d}x.$$ 
\par
Following the argument in Yang \cite[p.4]{yang2022extreme}, we define
$$
\mcg(s) \coloneqq \mcg_{\mbk,\ell}(s) = 1 + \frac{a_{\mbk}(2)}{2^s} + (-1)^{\ell} \zeta_{\mbk}^{(\ell)}(s).
$$
According to the above definition, it can be written as the following Dirichlet series for $\Re(s) > 1$
\begin{equation}
    \label{Gdefinition}
    \mcg(s) = \sum_{n = 1}^{\infty} \frac{a(n)}{n^s},
\end{equation}
which is absolutely convergent. Clearly, $a(1)=1, a(2) = a_\mbk(2)+a_\mbk(2) (\log 2)^{\ell},$ and for $n \ge3,a(n)=a_\mbk(n)(\log n)^{\ell}.$ The term $1 + a_{\mbk}(2) 2^{-s}$ added in the above definition ensures that the inequality $a(n) \ge a_\mbk(n)$ holds for all $n \ge 1$. Moreover, we have the trivial upper bound
\begin{equation}
    \label{tirvialupper}
    \mcg(s) \ll 1 + |\zeta_\mbk^{(\ell)}(s)|.
\end{equation}
\par
We now outline the structure of this paper. Section \ref{preliminarylemmas} presents several auxiliary lemmas that will be used throughout the proofs. In Section \ref{proofofthm1.1}, we construct an appropriate resonator and prove Theorem \ref{thm1.1}. The proof of Theorem \ref{thm1.2} will be given in Section \ref{proofofthm1.2}, which is preceded by the construction of another resonator. In the final Section \ref{supplementaryconclusions}, we present additional results that further refine our study of extreme values of derivatives of $\zeta_\mbk(s)$ within the critical strip.

\section{Preliminary Lemmas}
\label{preliminarylemmas}
In this section, we present several lemmas for late use.
The first lemma provides a growth estimate for derivatives of $\zeta_\mbk(s)$ in the critical strip.
\begin{lem}
    \label{approximatededekind}
    Fix $\delta \in (0, 1/2]$. Then uniformly for all $|t| \ge 1$ and $-\delta \le \sigma \le 1+\delta$, we have
    $$
    \zeta_\mbk^{(\ell)} ( \sigma + it) \ll d^{\phi(d)}|t| ^{\phi(d)(1-\sigma+\delta)/2}.
    $$
\end{lem}
\begin{proof}
    This estimate follows directly from \cite[Theorem 4]{rademacher1959phragmen} together with Cauchy's integral formula.
\end{proof}
\begin{xrem}
    A sharper bound for derivatives of $\zeta_\mbk(s)$ on the critical line was given in \cite{heath1988growth}. According to the Phragm\'en-Lindel\"of theorem, the exponent can be optimized. However, Lemma \ref{approximatededekind} is already sufficient for our purpose.
\end{xrem}
Furthermore, in order to establish a connection with GCD sums, we need to use a double-version convolution formula, following the approach of de la Bretèche and Tenenbaum \cite[Lemma 5.3]{tenen2019galsum}. It is similar to that in \cite[Lemma 3]{yang2022extreme}. For the sake of completeness, we provide a proof.

\begin{lem}
    	\label{convolutionformula}
      Suppose $\sigma \in [1/2, 1)$ and let $s = x + iy$. Let $K(s)$ be holomorphic in the horizontal strip $\sigma-2 \le y \le 0,$ satisfying the growth condition
     $$
     \max_{\sigma-2\le y \le 0}|K(x+iy)|=O\Big(\frac{1}{|x|^{2 \phi(d)} + 1}\Big), \quad |x| \to \infty.
     $$
     Then for all $t \neq 0,$ we have
     $$
     \int_{-\infty}^{\infty} \mcg(\sigma + it + iy)\mcg(\sigma - it + iy)K(y)\mathrm{d}y = \sum_{m,n\ge1} \frac{\widehat{K}(\log nm) a(n)a(m)}{n^{\sigma + it} m^{\sigma - it}}  - (\tau^{+} + \tau^{-}),
     $$
     where $\tau^{+}$ and $\tau^{-}$ are given by
     $$
     \tau^\pm = 2\pi \ell!\sum_{\substack{m + n = \ell \\ m,n\ge 0}} \frac{1}{m!n!} \Big. \Big( \frac{\mathrm{d}}{\mathrm{d}s} \Big)^m \mcg(s \pm it) \Big|_{s = 1 \pm it} \cdot \Big. \Big( \frac{\mathrm{d}}{\mathrm{d}s} \Big)^n K(i\sigma - is) \Big|_{s = 1 \pm it}.
     $$
\end{lem}
\begin{proof}
    Let $g(s) \coloneqq  \mcg(s + it )\mcg(s - it)K(i\sigma - is)$. The only poles in $\sigma \le {\textrm{Re}}(s) \le 2$ occur at $1 \pm it$.
    \par
    Integrating $g(s)$ along the rectangle with vertices $\sigma \pm iY$ and $2 \pm iY$, and applying the residue theorem, yields
    \begin{equation}
        \label{fourintegrals}
            I = I_1 + I_2 + I_3 + I_4 = \tau^{+} + \tau^{-},
    \end{equation}
    where $I_1$ and $I_3$ are the integrals of $g(s)$ from $\sigma - iY$ to $2 - iY$ and from $2 + iY$ to $\sigma + iY$ respectively, along the horizontal edges of the rectangle. Meanwhile, $I_2$ and $I_4$ are the integrals of $g(s)$ from $2 - iY$ to $2 + iY$ and from $\sigma + iY$ to $\sigma - iY$ respectively, corresponding to the vertical direction. We estimate these four integrals one by one in the subsequent calculations.
    \par
    It is clear that 
    $$
    \lim_{Y \to \infty} I_4 = -i \int_{-\infty}^{\infty} \mcg(\sigma + it + iy)\mcg(\sigma - it + iy)K(y)\mathrm{d}y,
    $$
    the infinite integral on the right-hand side is absolutely convergent due to the rapid decay of $K$. Moreover, by Cauchy's theorem, we conclude that
    $$
    \lim_{Y \to \infty} I_2 = i \sum_{n = 1}^{\infty} \sum_{m = 1}^{\infty}  \frac{\widehat{K}(\log nm)a(n)a(m)}{n^{\sigma + it} \cdot m^{\sigma - it}}.
    $$
    \par
    Applying \eqref{tirvialupper}, Lemma \ref{approximatededekind} and the rapid decay of $K$, we find that both $I_1$ and $I_3$ are bounded by
    $$
   Y^{-2 \phi(d)}\int_{\sigma}^{2} (1 + Y^{(1 - v + 3\delta)\phi(d)}) \mathrm{d}v \ll Y^{-2 \phi(d)} \Big( 1 + \frac{Y^{(1 + 3\delta - \sigma)\phi(d)} - Y^{2\delta\phi(d)}}{\log Y} \Big).
    $$
    Setting $\delta = 1/12$ ensures that both $I_1$ and $I_3$ tend to $0$ as $Y \to \infty$. Substituting this into \eqref{fourintegrals} completes the proof. 
\end{proof}
Similarly, we have the following single-version convolution formula.
\begin{lem}
    	\label{singleconvolutionformula}
     Suppose $\sigma \in [1/2, 1)$ and let $s = x + iy$. Let $K(s)$ be holomorphic in the horizontal strip $\sigma-2 \le y \le 0$, satisfying the growth condition
     $$
     \max_{\sigma-2\le y \le 0}|K(x+iy)|= O \Big(\frac{1}{|x|^{2 \phi(d)} + 1}\Big), \quad |x|\to \infty.
     $$
     Then for all $t \neq 0,$ we have
     $$
     \int_{-\infty}^{\infty} \mcg(\sigma + it + iy)K(y)\mathrm{d}y = \sum_{n = 1}^{\infty} \frac{\widehat{K}(\log n) a(n)}{n^{\sigma + it} }  -  2\pi \tau,
     $$
     where $\tau$ is given by
       $$
       \tau = \Big. \Big( \frac{\mathrm{d}}{\mathrm{d}s} \Big)^{\ell}K(i\sigma - i s)  \Big|_{s = 1 - it}.
       $$
\end{lem}
\begin{proof}
     The proof follows similarly from Lemma \ref{convolutionformula} by considering $g_1(s) \coloneqq \mcg(s+it)K(i\sigma-is)$. Therefore, the details are omitted for brevity.
\end{proof}
We need an appropriate smooth kernel function $K$ when applying Lemmas \ref{convolutionformula} and \ref{singleconvolutionformula}. Following the ideas in \cite{bondarenko2023dichotomy,bondarenko2018argument}, we set 
\begin{equation}
\label{Ketaudefinition}
    K_{\eta}(u) \coloneqq \frac{\sin^{2\eta}((\varepsilon\eta\inv \log T)u)}{(\varepsilon\eta\inv\log T)^{2\eta - 1}u^{2\eta}},
\end{equation}
which is the same choice as in \cite{bondarenko2023dichotomy}. Here $\varepsilon$ is small and $\eta \in \mbn$ will be chosen later. We present some properties of $\widehat{K_\eta}$, which will be used frequently in the subsequent proof. 

\begin{lem}[\cite{bondarenko2023dichotomy}, Lemma 5]
   \label{Ketauproperties}
Let $K_{\eta}(u)$ be defined as in \eqref{Ketaudefinition}. Then $\widehat{K_{\eta}}(v)$ is a real and even function supported on $|v|\le 2\varepsilon\log T.$ 
It is decreasing on $[0,\infty)$ and satisfies $0\le\widehat{K_{\eta}}(v)\le\widehat{K_{\eta}}(0)$. Moreover, we have
\begin{equation}
\label{Ketahatderivatives}
\Big|\frac{\mathrm{d}}{\mathrm{d}v}\widehat{K_{\eta}}(v)\Big|\le\frac{\widehat{K_{\eta-1}}(0)}{\varepsilon\eta\inv\log T}.
\end{equation}
Furthermore, for large $\eta$, the following equivalence holds:
\begin{equation}
\label{Ketahatequivalence}
\widehat{K_{\eta}}(0)\sim\sqrt{\frac{3\pi}{\eta}}.
\end{equation}
\end{lem}

\section{Proof of Theorem \ref{thm1.1}}
\label{proofofthm1.1}
In this section, we apply the auxiliary lemmas from Section \ref{preliminarylemmas} to complete the proof of Theorem \ref{thm1.1}. We first establish lower bounds for GCD sums weighted by the coefficients $a_\mbk(n)$.
\subsection{Lower bounds for Weighted GCD sums}
We extend the construction introduced by de la Bretèche and Tenenbaum \cite{tenen2019galsum}. Define
\begin{equation}
\label{agcddefine}
S_{1/2}(\mcm,a_\mbk) \coloneqq \sum_{m,n\in \mcm} a_\mbk\Big(\frac{m}{(m,n)}\Big)a_\mbk\Big(\frac{n}{(m,n)}\Big)\sqrt{\frac{(m,n)}{[m,n]}}. 
\end{equation}
\par
Let $\alpha \in (1,e)$ and $\delta \in (0,1)$ be some parameters, and let $N = \lfloor T^{\lambda} \rfloor$ be the integer not larger than $T$ with $\lambda$ to be chosen later. Let $k$ be an integer satisfying $1\le k \le (\log_{2}N)^{\delta},$ then define $P_k$ as the set of primes $p$ that satisfy 
$$
 \phi(d)\alpha^{k}\log N \log_{2}N < p \le  \phi(d)\alpha^{k + 1}\log N \log_{2}N \quad \text{and} \quad p\equiv 1  \pmod{d}.
$$
Since $d \ll (\log_2 T)^A $, the Siegel-Walfisz theorem implies 
$$
 |P_{k}| \sim \alpha^{k}(\alpha - 1)\log N \bigg(1 + O\Big(\frac{k+\log_{3}N}{\log_{2}N}\Big)\bigg).
$$
Let $\beta > 1$ satisfy $\beta \delta \log \alpha < 1$. We then define 
$$
W_k \coloneqq 2\Big\lfloor\frac{\beta\log N}{2k^{2}\log_{3}N}\Big\rfloor.
$$ 
Furthermore, we let $W \coloneqq \prod_{p\in P_{k}}p$, and consider the following set
$$
\mcm_k\coloneqq \Big\{m:m = \frac{\ell}{q}W, \ \omega(\ell)\le \frac{1}{2}W_{k}, \omega(q) \le \frac{1}{2}W_{k}, \ \ell q\mid W\Big\},
$$
where $\omega(n)$ denotes the number of distinct prime factors of $n$. Next, we consider
\begin{equation}
    \label{Mdefinition}
    \mcm\coloneqq\bigg\{m = \prod_{1\leq k\leq (\log_{2}N)^{\delta}}m_{k}: \ m_{k}\in\mcm_k \, , \, 1\leq k\leq (\log_{2}N)^{\delta}\bigg\}.
\end{equation}
\par
According to \cite[Lemma 2.2]{tenen2019galsum}, we have $|\mcm|\le N$. We establish the following lower bound for $S_{1/2}(\mcm,a_\mbk)$, and the proof will be given in the Appendix.
\begin{prop}
    \label{agcdlowerbound}
    Define $S_{1/2}(\mcm,a_\mbk)$ as in \eqref{agcddefine}, where $\mcm$ is given in \eqref{Mdefinition}. Then for $d \ll (\log_2 T)^A$, we have
    $$\frac{S_{1/2}(\mcm, a_\mbk)}{|\mcm|} \ge \exp\bigg(\Big(2\sqrt{\phi(d)} + o(1)\Big)\sqrt{\frac{\log N \log_3 N}{\log_2 N}} \bigg)$$
    as $N \to \infty$.
\end{prop}
Fonga attempted to extend the results of de la Bretèche and Tenenbaum in \cite[Theorem 1.1]{fonga2023extreme}. We argue that the range of $k$ he chose is inappropriate. Accordingly, we have revised it and supplemented the proof in the Appendix.

\begin{xrem}
    For $1/2 < \sigma <1$, one can define analogously:
    \begin{equation}
            \label{SsigmamathcalMa}
        S_{\sigma}(\mcm,a_\mbk) \coloneqq \sum_{m,n\in \mcm} a_\mbk\Big(\frac{m}{(m,n)}\Big)a_\mbk\Big(\frac{n}{(m,n)}\Big)\Big(\frac{(m,n)}{[m,n]}\Big)^{\sigma}.
    \end{equation}
   If a lower bound similar to that in Proposition \ref{agcdlowerbound} can be obtained for $S_{\sigma}(\mcm,a_\mbk)$, the result of Theorem \ref{thm1.2} can be optimized through a process nearly analogous to that of Theorem \ref{thm1.1}.
\end{xrem}

\subsection{Constructing the resonator}
\label{resonator1}
Following the ideas from \cite{aistleitner2016lower,bondarenko2017large,tenen2019galsum,yang2022extreme}, we construct the resonator $R(t)$ to apply the resonance method. For each integers $j \ge 0$, define
$$
\mcn_j\coloneqq\bigg[\Big(1 + \frac{\log T}{T}\Big)^{j}, \Big(1 + \frac{\log T}{T}\Big)^{j + 1}\bigg)\cap\mcm .
$$
Here $T$ is a parameter, and $\mcm$ is defined as \eqref{Mdefinition}. Using the notations in \cite{tenen2019galsum}, we let $h_j \coloneqq \min\mcn_j$ if $\mcn_j \neq \emptyset$. Define $\mch$ as the set of all $h_j$ and we consider the function $r$ on $\mch$ defined by
$$
r(h_{j})=\sqrt{\sum_{m\in\mcn_j}1}.
$$
We then define the resonator as follows:
$$
R(t)=\sum_{h\in\mathcal{H}}r(h)h^{-it}.
$$
By the Cauchy-Schwarz inequality, we have
$$
|R(t)|^2 \le R(0)^{2}\le N\sum_{h\in\mch}r(h)^{2}\le N|\mcm|\le N^{2}.
$$
\par
Set $\Phi(y)\coloneqq e^{-y^{2}/2}$ as in \cite{bondarenko2017large}. Plainly, the Fourier transform $\widehat{\Phi}$ satisfies $\widehat{\Phi}(\xi)=\sqrt{2\pi}\Phi(\xi)$. According to \cite[Lemma 5]{bondarenko2018argument}, we obtain
\begin{equation}
\label{trivialboundrphi}
\int_\mbr |R(t)|^{2}\Phi\Big(\frac{t\log T}{T}\Big)\mathrm{d}t\ll\frac{|\mcm|T}{\log T}.
\end{equation} 

\subsection{The proof of Theorem \ref{thm1.1}}
\label{proofthm1.1process}
We now complete the proof of Theorem \ref{thm1.1}. Let $\eta = 2\phi(d)$ as in \cite{bondarenko2023dichotomy} and $N = \lfloor T^\lambda \rfloor$. For fixed $\varepsilon >0$, we choose $\lambda$ such that $\lambda +5\varepsilon <1$. We start with the following integral
$$
J(T)\coloneqq\int_{|t|\ge1}\int_\mbr\mcg\Big(\frac{1}{2} + it + iu\Big)\mcg\Big(\frac{1}{2} - it + i u\Big)K_{\eta}(u)|R(t)|^{2}\Phi\Big(\frac{t\log T}{T}\Big)\mathrm{d}u\mathrm{d}t.
$$
We divide $J(T)$ into three parts as follows
\begin{align}
\label{J=J123}
J(T) & = \int_{1\le |t| \le T \log T}\int_{|u| \le T^{1/2}} + \int_{1\le |t| \le T \log T}\int_{|u| \ge T^{1/2}} + \int_{|t| \ge T \log T}\int_{\mbr} \nonumber \\
& \eqqcolon J_1 (T)+ J_2 (T)+ J_3 (T).
\end{align}
For brevity, the integrands on the right-hand side are omitted, without causing any ambiguity. We will show that the integral $J_1(T)$ provides the main term for $J(T)$.
\par
Using Lemma \ref{approximatededekind}, we establish the following bound 
\begin{equation}
        \label{GKellupperbound}
\Big|\mcg\Big(\frac{1}{2}\pm it + iu\Big)\Big|\ll 1 + \Big|\zeta_\mbk^{(\ell)}\Big(\frac{1}{2} \pm it + iu\Big)\Big|\ll (1 + |t| + |u|)^{\frac{3}{10}\phi(d)}.
\end{equation}
Combining with $K_{\eta}(u) \ll u^{-4\phi(d)},$ it is clear that
\begin{align*}
    J_2 (T) &\,\ll d^{\phi(d)} \int_{1\le |t| \le T \log T}\int_{|u| \ge T^{1/2}} (1 + |t| + |u|)^{3\phi(d)/5} K_{\eta}(u)|R(t)|^{2}\Phi\Big(\frac{t\log T}{T}\Big)\mathrm{d}u\mathrm{d}t \\
           &\,\ll d^{\phi(d)} ( T \log T)^{3\phi(d)/5}\int_{|u| \ge T^{1/2}} |u| ^{3\phi(d)/5-4\phi(d)}\mathrm{d}u \int_\mbr|R(t)|^{2}\Phi\Big(\frac{t\log T}{T}\Big)\mathrm{d}t.
\end{align*}
Then, by \eqref{trivialboundrphi} and $d \ll (\log_2 T)^A,$
\begin{equation}
\label{J2upperbound}
    J_2(T)  \ll d^{\phi(d)} |\mcm|T^{-7\phi(d)/20+\varepsilon}=o(|\mcm|T).
\end{equation}
Moreover,
$$
J_3 (T) \ll d^{\phi(d)} \int_{|t| \ge T \log T} |t|^{3\phi(d)/5}|R(t)|^{2}\Phi\Big(\frac{t\log T}{T}\Big)\mathrm{d}t \int_\mbr |u|^{3\phi(d)/5} K_{\eta}(u)\mathrm{d}u.
$$
The rapid decay of $\Phi$ and \eqref{GKellupperbound} implies that
\begin{equation}
    \label{J3upperbound}
    J_3 \ll d^{\phi(d)}|\mcm| T^{3\phi(d)/5+1+\varepsilon} \Phi \Big(\frac{(\log T)^2}{\sqrt{2}}\Big)=o(|\mcm|T).
\end{equation}
Substituting \eqref{J2upperbound}, \eqref{J3upperbound} into \eqref{J=J123}, we have
\begin{equation}
    \label{JJ1+o}
    J(T) = J_1 (T) + o(|\mcm|T) .
\end{equation}
\par
Define
$$
\mcz \coloneqq \mcz_\mbk(T)= \max_{-T \log T - T^{1/2}\le t\le T \log T+T^{1/2}}\Big|\mcg\Big(\frac{1}{2}+it\Big)\Big|.
$$
Combining with \eqref{trivialboundrphi}, we immediately derive the following lower bound
\begin{equation}
    \label{J1upperbound}
    J_1(T) \ll \widehat{K_{\eta}}(0) \mcz^2 \frac{|\mcm|T}{\log T}.
\end{equation}
Thus, 
\begin{equation}
    \label{Jupperbound}
    J(T) \ll \widehat{K_{\eta}}(0) \mcz^2 \frac{|\mcm|T}{\log T} +  o (|\mcm|T).
\end{equation}
\par
Next, we seek an effective lower bound for $J(T)$. Our approach establishes a connection between $J(T)$ and GCD sums via the double-version convolution formula. We introduce the following three integrals:
\begin{align*}
&J_0(T)\coloneqq \int_{|t|\ge 1}E(t)|R(t)|^{2}\Phi\Big(\frac{t\log T}{T}\Big)\mathrm{d}t, \\
&\mce^+\coloneqq -\int_{|t|\ge 1}\tau^{+}\cdot|R(t)|^{2}\Phi\Big( \frac{t\log T}{T}\Big)\mathrm{d}t, \\
&\mce^-\coloneqq -\int_{|t|\ge 1}\tau^{-}\cdot|R(t)|^{2}\Phi\Big( \frac{t\log T}{T}\Big)\mathrm{d}t,
\end{align*}
where the series $E(t)$ is defined as 
\begin{equation}
    \label{EKtseriesdefinition}
E(t)\coloneqq \sum_{m,n\ge1}\frac{\widehat{K_{\eta}}(\log nm)}{\sqrt{nm}(n/m)^{it}}a(n)a(m).
\end{equation}
Lemma \ref{convolutionformula} implies that $J(T)=J_0(T)+\mce^+ +\mce^-.$ As in \cite[pp. 14-15]{yang2022note}, we get
$$
|\mce^+|\ll|\mcm| \frac{T^{\lambda}}{\log T} \quad \text{and} \quad |\mce^-|\ll|\mcm| \frac{T^{\lambda}}{\log T}.
$$
Subsequently, we use the Fourier transform on the real axis and therefore focus on the region $|t| \le 1$ now. Lemma \ref{Ketauproperties} shows that $\widehat{K_{\eta}}(\log nm)=0$ when $mn\ge T^{2\varepsilon}$. Thus, we obtain
$$
\int_\mbr E(t)|R(t)|^2\Phi\Big(\frac{t\log T}{T} \Big)\mathrm{d}t \ll \widehat{K_\eta}(0)R(0)^2\sum_{m,n\ge1}\frac{a(n)a(m)}{\sqrt{mn}}.
$$
Since $d \ll (\log_2 T)^A,$ we have $a(n)n^{-1/2}\ll (\log_2 T)^A.$ Thus, the integral on $|t| \le1$ can be bounded by
$$
\widehat{K_{\eta}}(0)R(0)^2(\log_2 T)^A \sum_{mn\le T^{2\varepsilon}}1 \ll \widehat{K_{\eta}}(0)|\mcm|T^{\lambda + 5\varepsilon}.
$$
We extend the integral to $\mbr$. Specifically, setting
$$
\widetilde{J}(T)\coloneqq\int_\mbr E(t)|R(t)|^{2}\Phi\Big(\frac{t\log T}{T}\Big)\mathrm{d}t,
$$
we have
$$
J_0(T) = \widetilde{J}(T) + O(|\mcm|T^{\lambda + 5\varepsilon}).
$$
Hence, combining with \eqref{Jupperbound}, we obtain 
\begin{equation}
    \label{Jtuitaupperbound}
    \widetilde{J}(T) \ll  \widehat{K_{\eta}}(0) \mcz^2 \frac{|\mcm|T}{\log T} +  o(|\mcm| T) + O(|\mcm|T^{\lambda + 5\varepsilon}).
\end{equation}
\par
We now derive the lower bound for $\widetilde{J}(T)$. Expanding $|R(t)|^2$ and $E(t)$ in $\widetilde{J}(T)$ yields that 
$$
\widetilde{J}(T)=\sqrt{2\pi}\frac{T}{\log T}\sum_{g,k\geq1}\frac{\widehat{K_{\eta}}(\log gk)a(g)a(k)}{\sqrt{gk}}\sum_{h,h^{\prime}\in\mch}r(h)r(h^{\prime})\Phi\Big(\frac{T}{\log T}\log\frac{hg}{h^{\prime}k}\Big).
$$
For the terms involving $\widehat{K_{\eta}}(u)$, we follow the approach in \cite[p. 10]{bondarenko2023dichotomy}, using \eqref{Ketahatderivatives} and \eqref{Ketahatequivalence}.
For the inner sum, it follows from the proof of \cite[Theorem 1.4]{tenen2019galsum} that
$$
\widetilde{J}(T) \gg \widehat{K_{\eta}}(0)\frac{T}{\log T}\sum_{\substack{m,n\in\mcm \\  mg = nk}}\sum_{1\le gk\le T^{\varepsilon/3\eta}}\frac{a_\mbk(g)a_\mbk(k)}{\sqrt{gk}}.
$$
\par
Setting 
$$
g = \frac{n}{(m,n)}, \quad k = \frac{m}{(m,n)}, \quad gk = \frac{[m,n]}{(m,n)} \le T^{\varepsilon/3\eta},
$$
we have 
$$
\widetilde{J}(T) \gg \widehat{K_{\eta}}(0)\frac{T}{\log T}\sum_{\substack{m,n\in\mcm\\ [m,n]/(m,n)\le T^{\varepsilon/3\eta}}}a_\mbk\Big(\frac{m}{(m,n)}\Big)a_\mbk\Big(\frac{n}{(m,n)}\Big)\sqrt{\frac{(m,n)}{[m,n]}}.
$$
Using Rankin's trick, we can establish a connection between $\widetilde{J}$ and weighted GCD sums as follows
$$
\widetilde{J}(T)\gg \widehat{K_{\eta}}(0)\frac{T}{\log T}(S_{1/2}(\mcm, a_\mbk) - T^{-\varepsilon/18\eta}S_{1/3}(\mcm,a_\mbk)),
$$
where $S_{1/2}(\mcm,a_\mbk)$ and $S_{1/3}(\mcm,a_\mbk)$ are given by \eqref{agcddefine} and \eqref{SsigmamathcalMa}, respectively.
\par
Combining with \eqref{Jtuitaupperbound}, we have
$$
\mcz^2  \gg \frac{S_{1/2}(\mcm, a_\mbk)}{|\mcm|} - \frac{S_{1/3}(\mcm, a_\mbk)}{|\mcm|T^{\varepsilon/18\eta}} + o(\log T).
$$
Following the idea from \cite{tenen2019galsum}, if we set $P^+ (m)$ as the largest prime factor of $m$ and $y_{\mcm} \coloneqq \max_{m \in \mcm} P^+(m),$ then we have 
$$\frac{S_{1/3}(\mcm, a_\mbk)}{|\mcm|} \ll \exp \Big(\phi(d) y_\mcm^{2/3}\Big).$$
Based on the construction of $P_k$ in Section \ref{preliminarylemmas}, we have $y_\mcm \ll (\log T)^{6/5}$, so that for sufficiently large $T$, 
$$
\frac{S_{1/3}(\mcm, a_\mbk)}{T^{\varepsilon/18\eta}|\mcm|} \ll \frac{\exp\Big(\phi(d) y_\mcm^{2/3}\Big)}{T^{\varepsilon/18\eta}} \to 0.
$$
Thus, by Proposition \ref{agcdlowerbound}, for all sufficiently large $T$, we have 
$$
\mcz^2 \gg \exp \bigg(\Big(2\phi(d)\sqrt{\lambda}+o(1)\Big)\sqrt{\frac{\log T\log_3 T}{\log_2 T}}\bigg).
$$
Trivially, 
\begin{align}
\label{maxtlog+t1/2}
    \max_{-T \log T - T^{1/2}\le t\le T \log T +T^{1/2}}&\,\Big|\zeta_\mbk^{(\ell)}\Big(\frac{1}{2}+it\Big)\Big| \gg \mcz + O(1) \nonumber \\
    &\, \gg \exp\bigg(\Big(\phi(d)\sqrt{\lambda}+o(1)\Big)\sqrt{\frac{\log T\log_3 T}{\log_2 T}}\bigg).
\end{align}
\par
We let $\lambda$ be sufficiently close to 1 such that $\lambda + 5\varepsilon <1$ for small $\varepsilon$. Finally, we set 
$$T^{\prime}=(1+o(1))\frac{T}{\log T}$$
and readjust the parameter $T=T^{\prime}$ following the idea from \cite{bondarenko2023dichotomy}. Variations in the logarithmic factor affect only the lower order terms in the exponent on the right-hand side of \eqref{maxtlog+t1/2}. Hence, this result also holds for $t \in [0,T]$. We complete the proof of Theorem \ref{thm1.1}.

\section{Proof of Theorem \ref{thm1.2}}
\label{proofofthm1.2}
In this section, we prove Theorem \ref{thm1.2} using the resonance method. As in Section \ref{proofofthm1.1}, we set $\eta = 2 \phi(d)$.
\subsection{Constructing the resonator}
\label{resonator2}
We construct the resonator $R(t)$ following a similar strategy as in \cite{chirre2019extreme}. Let $\gamma \in (0,1)$ be a parameter to be chosen later. Define $\mcp_d$ as the set of prime numbers $p$ such that
$$
e \phi(d) \log N \log_{2}N < p \le \phi(d)\log N\exp((\log_2 N)^{\gamma})\log_{2}N \quad \text{and} \quad p\equiv 1 \pmod{d}.
$$
For any $1 \le c_d \le \sqrt{\phi(d)/(e-1)}$, define the multiplicative function $f(n)$ supported on the set of square-free numbers, by setting
$$
f(p) = c_d \bigg( \frac{(\log N)^{1 - \sigma} (\log_2 N)^{\sigma}}{(\log_3 N)^{1 - \sigma}} \bigg) \frac{1}{p^{\sigma} (\log p - \log_2 N - \log_3 N - \log \phi(d))}
$$
as $p \in \mcp_d$, and $f(p) = 0$ otherwise.
\par
For each integer $k \in \{1, \cdots, \lfloor(\log_2 N)^{\gamma}\rfloor\}$, define $\mcp_{k,d}$ as the set of $p$ satisfying
$$
e^k \phi(d) \log N \log_{2}N < p \le e^{k+1}\phi(d)\log N\log_{2}N \quad \text{and} \quad p\equiv 1 \pmod{d},
$$
which is the subset of $\mcp_d$. Fixing $1 < \alpha < 1/\gamma$, define 
$$
\mcm_{k,d} \coloneqq \bigg\{ n \in \mathrm{supp}(f) : n \text{ has at least } \Delta_k \coloneqq \frac{\alpha (\log N)^{2 - 2\sigma}}{k^2 (\log_3 N)^{2 - 2\sigma}} \text{ prime divisors in } \mathcal{P}_{k,d} \bigg\}.
$$ 
Let $\mcm^{\prime}_{k,d}$ consist of those integers $\mcm_{k,d}$ whose  prime factors all lie in $\mcp_{k,d}$. Subsequently, set
$$
\mcm_d \coloneqq \mathrm{supp}(f) \setminus \bigcup_{k = 1}^{\lfloor(\log_2 N)^{\gamma}\rfloor} \mcm_{k,d}. 
$$
\par
From the above definitions, we can deduce several useful properties. First, $\mcm_d$ is divisor closed. Specifically, if some $m^{\prime} \mid m \in \mcm_d,$ then $m^{\prime} \in \mcm_d$. Adapting the argument of \cite[Lemma 8]{chirre2019extreme}, we obtain that $|\mcm_d| \le N$. Furthermore, the Siegel-Walfisz theorem gives
\begin{equation}
\label{Pkdcardinality}
  |\mcp_{k,d}| \le (1 + o(1)) e^{k+1} \log N.  
\end{equation}
\par
Now we construct the resonator $R(t)$. Define $\mcj$ as the set of integers $j$ such that
$$
\bigg[\Big(1 + \frac{\log T}{T}\Big)^j, \Big(1 + \frac{\log T}{T}\Big)^{j + 1}\bigg] \bigcap \mcm_d \neq \emptyset,
$$ 
and let $m_j$ be the minimum of $[(1 + \log T/T)^j, (1 + \log T/T)^{j + 1}] \cap \mcm_d $ for all $j \in \mcj.$ Setting $\mcm^{\prime}_d = \{m_j : j \in \mcj \}$, define
$$
r(m_j) \coloneqq \bigg( \sum_{n \in \mcm_d, (1 - \log T/T)^{j - 1} \le n \le (1 + \log T/T)^{j + 2}} f(n)^2 \bigg)^{1/2}.
$$
\par
Finally, define 
$$
R(t) \coloneqq \sum_{m \in \mcm^{\prime}_d} r(m) m^{-it}.
$$
By the Cauchy-Schwarz inequality,
$$
|R(t)|^2 \le R(0)^2 \ll N \sum_{n \in \mbn} f(n)^2.
$$
Moreover, taking $\Phi$ to be the same as in Section \ref{proofofthm1.1}, and using \cite[Lemma 5]{bondarenko2018argument}, we obtain
$$
\int_\mbr |R(t)|^{2}\Phi\Big(\frac{t\log T}{T}\Big)\mathrm{d}t \ll \frac{T}{\log T}\sum_{n \in \mbn} f(n)^2.
$$

\subsection{Supporting propositions}
\label{supportingpropositions}
Define
\begin{align*}
    \mca_d &\coloneqq \mca_{d,N}= \frac{1}{\sum_{n\in \mbn} f(n)^2} \sum_{n\in\mbn} \frac{f(n)}{n^{\sigma}} \sum_{q\mid n} a_\mbk\Big(\frac{n}{q}\Big)f(q)q^{\sigma} \\
    & = \prod_{p\in \mcp_d} \frac{1 + f(p)^2 + \phi(d)f(p)p^{-\sigma}}{1 + f(p)^2}.
\end{align*}
Here, $a_\mbk(n)$ are the coefficients of the Dirichlet series of $\zeta_\mbk(s)$, and $f$ is the multiplicative function introduced in Section \ref{resonator2}. In this subsection, we present 
some auxiliary propositions that provide sharp lower bound for $\mca_d$ and show that the error terms caused by truncation or the exclusion of $\mcm_d$ are negligible. These results play a central role in establishing the lower bound in Theorem \ref{thm1.2}.
\begin{prop}
    \label{prop41}
    Let $d \ll (\log_2 N)^A$ and $c_d \le \sqrt{\phi(d)/(e-1)}$. Then for $0 < \delta <1$ as $N \to +\infty$, we have
    $$\mca_d \geq \exp\bigg((\delta\gamma c_d+ o(1))\frac{(\log N)^{1-\sigma}(\log_3 N)^{\sigma}}{(\log_2 N)^{\sigma}}\bigg).$$
\end{prop}
\begin{proof}
    For $c_d \le \sqrt{\phi(d)/(e-1)}$, we have $f(p)=o(1)$ for all $p \in \mcp_d$. It follows that
    $$
    \mca_d = \exp \Big((1+o(1)) \sum_{p \in \mcp_d} \frac{f(p)}{p^{\sigma}}\Big).
    $$
    By the definition of $f(p)$, we have
      \begin{align*}
      \sum_{p \in \mcp_d} \frac{f(p)}{p^{\sigma}} 
      & =c_d \frac{(\log N)^{1 - \sigma} (\log_2  N)^{\sigma}}{(\log_3 N)^{1 - \sigma}} \sum_{p \in \mcp_d} \frac{1}{p^{2\sigma}(\log p-\log _2 N-\log _3 N - \log \phi(d))}  \\
      & \ge c_d \frac{(\log N)^{1 - \sigma} (\log_2 N)^{\sigma}}{(\log_3 N)^{1 - \sigma}} \sum_{k=1}^{\lfloor(\log_2 N)^{\gamma}\rfloor} \sum_{p \in \mcp_{k,d}} \frac{1}{(k+1)p^{2\sigma}}.
      \end{align*}
      Notice that \cite[Lemma 9]{chirre2019extreme}, we get the following lower bound
    $$
    \sum_{p \in \mcp_{k,d}} \frac{1}{p^{2\sigma}} > (\delta + o(1))(\log_2 N)^{ -2\sigma},
    $$
    where $0< \delta <1$. Thus,
    $$
    \sum_{p \in \mcp_d} \frac{f(p)}{p^{\sigma}} \ge (\delta\gamma c_d + o(1))\frac{(\log N)^{1-\sigma} (\log_3 N)^{\sigma}}{(\log_2 N)^{\sigma}}.
    $$
    This completes the proof.
\end{proof}

\begin{prop}
    \label{prop42}
    Let $d \ll (\log_2 N)^A$ and $c_d \le \sqrt{\phi(d)/(e-1)}$. Then     as $N \to +\infty$, we have
    $$
    \frac{1}{\sum_{n \in \mbn} f(n)^2} \sum_{n \in \mbn, n \notin \mcm_d} \frac{f(n)}{n^{\sigma}} \sum_{q \mid n} a_\mbk \Big( \frac{n}{q}\Big)f(q)q^{\sigma}=o(\mca_d).
    $$
\end{prop}
\begin{proof}
By the definition of $\mcm_d,$ we obtain
\begin{align*}
\mcf &\coloneqq \frac{1}{\mca_d \sum_{n\in\mbn} f(n)^2} \sum_{n\in\mbn, n\notin\mcm_d} \frac{f(n)}{n^{\sigma}} \sum_{q\mid n}a_\mbk \Big( \frac{n}{q}\Big)f(q)q^{\sigma} \\
&\le \frac{1}{\mca_d \sum_{n\in\mbn} f(n)^2} \sum_{k = 1}^{\lfloor(\log_2 N)^{\gamma}\rfloor} \sum_{n\in\mcm_{k,d}} \frac{f(n)}{n^{\sigma}} \sum_{d\mid n} a_\mbk \Big( \frac{n}{q}\Big)f(q)q^{\sigma}.
\end{align*}
From the construction of $f(n)$ and $\mcm_{k,d}^{\prime}$, it follows that $\mcf$ can be bounded by
\begin{align*}
&\sum_{k = 1}^{\lfloor(\log_2 N)^{\gamma}\rfloor} \frac{1}{\prod_{p\in \mcp_{k,d}} (1 + f(p)^2 + \phi(d)f(p)p^{-\sigma})} \sum_{n\in\mcm^{\prime}_{k,d}} \frac{f(n)}{n^{\sigma}} \sum_{q\mid n} a_\mbk \Big( \frac{n}{q}\Big)f(q)q^{\sigma} \\
&\le \sum_{k = 1}^{\lfloor(\log_2 N)^{\gamma}\rfloor} \frac{1}{\prod_{p\in \mcp_{k,d}} (1 + f(p)^2)} \sum_{n\in\mcm^{\prime}_{k,d}} f(n)^2 \prod_{p\in \mcp_{k,d}} \Big(1 + \frac{\phi(d)}{f(p)p^{\sigma}}\Big) \eqqcolon \sum_{k = 1}^{\lfloor(\log_2 N)^{\gamma}\rfloor}  \mcf^{\star}.
\end{align*}
It is clear from \eqref{Pkdcardinality} that 
\begin{align*}
\prod_{p\in \mcp_{k,d}} \Big(1 + \frac{\phi(d)}{f(p)p^{\sigma}}\Big)& \le \prod_{p\in \mcp_{k,d}} \bigg(1 + \frac{\phi(d)}{c_d}(k+1)\frac{(\log_3 N)^{1-\sigma}}{(\log N)^{1-\sigma}(\log_2 N)^{\sigma}}\bigg)  \\
&\le \exp\bigg(\frac{\phi(d)}{c_d}(k + 1)e^{k + 1}\frac{(\log N)^{\sigma}(\log_3 N)^{1-\sigma}}{(\log_2 N)^{\sigma}}\bigg).
\end{align*}
Therefore, combining with the fact that $d \ll ( \log_2 N)^A$ and $k \le (\log_2 N)^{\gamma}$, we have 
$$
\prod_{p\in \mcp_{k,d}} \Big(1 + \frac{\phi(d)}{f(p)p^{\sigma}}\Big) =\exp\bigg(o\bigg(\frac{(\log N)^{2-2\sigma}}{k^2 (\log_3 N)^{2-2\sigma}}\bigg)\bigg).
$$
\par
On the other hand, recalling the definition of $\mcm^{\prime}_{k,d},$ we utilize the classical method in \cite[p. 8]{bondarenko2017large}. Since $f(n)$ is multiplicative, we obtain
$$
\sum_{n\in\mcm^{\prime}_{k,d}} f(n)^2 \le b^{-\Delta_k} \prod_{p\in \mcp_{k,d}} (1 + b f(p)^2)
$$
whenever $b>1$. Thus, 
$$
\frac{1}{\prod_{p\in \mcp_{k,d}} (1 + f(p)^2)} \sum_{n\in\mcm^{\prime}_{k,d}} f(n)^2 \le b^{-\Delta_k} \exp\Big( (b - 1)\sum_{p\in \mcm^{\prime}_{k,d}}f(p)^2\Big).
$$
For the sum in the exponent, the Siegel-Walfisz theorem yields 
\begin{align*}
\sum_{p\in \mcm^{\prime}_{k,d}} f(p)^2 &= c_d^2\frac{(\log N)^{2-2\sigma}(\log_2 N)^{2\sigma}}{(\log_3 N)^{2-2\sigma}} \sum_{p\in \mcp_{k,d}} \frac{1}{p^{2\sigma}(\log p - \log_2 N - \log_3 N - \log \phi(d))^2} \\
&\le \frac{c_d^2}{\phi(d)} (1 + o(1))\frac{(\log N)^{2-2\sigma}(\log_2 N)^{2\sigma}}{(\log_3 N)^{2-2\sigma}} \int_{e^k \phi(d) \log N \log_2 N}^{e^{k+1} \phi(d) \log N \log_2 N} \frac{\mathrm{d}x}{k^2 x^{2\sigma}\log x}\\
&\le \frac{c_d^2 (e-1)}{\phi(d)}(1 + o(1))\frac{(\log N)^{2-2\sigma}}{k^2 (\log_3 N)^{2-2\sigma}}.
\end{align*}
Combining all the above estimates, $\mcf^\star$ is bounded by
$$
\exp\bigg(\bigg(\frac{c_d^2 (e-1)}{\phi(d)}(b - 1) - \alpha\log b + o(1)\bigg)\frac{(\log N)^{2-2\sigma}}{k^2 (\log_3 N)^{2-2\sigma}}\bigg).
$$
Since $c_d \le \sqrt{\phi(d)/(e-1)}$ and $\alpha>1$, it suffices to choose $b$ sufficiently close to $1$ so that the exponent becomes negative. Hence, the proof is complete.
\end{proof}

\begin{prop}
    \label{prop43}
    Let $d \ll (\log_2 N)^A$ and $c_d \le \sqrt{\phi(d)/(e-1)}$. Then     as $N \to +\infty$, we have 
    $$
    \frac{1}{\sum_{n \in \mbn} f(n)^2} \sum_{n \in \mcm_d} \frac{f(n)}{n^{\sigma}} \sum_{q \mid n, q \leq n/N^{\varepsilon/3\eta}} a_\mbk \Big( \frac{n}{q}\Big)f(q)q^{\sigma}=o(\mca_d).
    $$
\end{prop}
\begin{proof}
We have
$$
\sum_{n \in \mcm_d} \frac{f(n)}{n^{\sigma}} \sum_{q \mid n, q \leq n/N^{\varepsilon/3\eta}} a_\mbk \Big( \frac{n}{q}\Big)f(q)q^{\sigma} = \sum_{n \in \mcm_d} f(n)^2 \sum_{q \mid n, q \geq N^{\varepsilon/3\eta}} \frac{a_\mbk(q)}{f(q)q^\sigma}.
$$
Furthermore, applying Rankin’s trick to the inner sum on the right-hand side yields
$$
\sum_{q \mid n, q \geq N^{\varepsilon/3\eta}} \frac{a_\mbk(q)}{f(q)q^\sigma} \le N^{-\varepsilon/12\eta}\sum_{q \mid n}  \frac{a_\mbk(q)}{f(q)q^{\sigma-1/4}} =  N^{-\varepsilon/12\eta} \prod_{p \mid n}\Big(1+ \frac{\phi(d)}{f(p)p^{\sigma-1/4}} \Big).
$$
Since $\eta = 2\phi(d)$, and because each $n$ in $\mcm_d$ has at most $ O((\log N / \log_3 N)^{2-2\sigma})$ prime factors, the proof follows by a process similar to that in \cite[Lemma 6]{bondarenko2023dichotomy}.
\end{proof}
\subsection{The proof of Theorem \ref{thm1.2}}
\label{proofthm1.2process}
In this subsection, the ideas and techniques in our preceding part are essentially the same as those used in Section \ref{proofthm1.1process}. Let $N=\lfloor T^\lambda \rfloor$, where $\lambda + 2 \varepsilon < 1$ holds for small $\varepsilon$. Define 
$$
W_1(T) \coloneqq \int_{1\le |t| \le T \log T}\int_{ |u| \le T^{1/2}}\mcg(\sigma + it + iu)K_{\eta}(u)|R(t)|^{2}\Phi\Big(\frac{t\log T}{T}\Big)\mathrm{d}u\mathrm{d}t
$$
and
$$
\mcz(\sigma) \coloneqq \mcz_\mbk(\sigma,T) =  \max_{-T \log T - T^{1/2}\le t\le T \log T+T^{1/2}}|\mcg(\sigma+it)|.
$$
Trivially,
$$
    W_1(T) \le
    \widehat{K_{\eta}}(0)\mcz(\sigma) \int_\mbr |R(t)|^{2}\Phi\Big(\frac{t\log T}{T}\Big)\mathrm{d}t  \ll \widehat{K_{\eta}}(0)\mcz(\sigma)\frac{T}{\log T}\sum_{n \in \mbn} f(n)^2.
$$
Then, we proceed by considering the tails of the integrals:
\begin{align*}
    W_2(T) &\, \coloneqq\int_{1\le |t| \le T \log T}\int_{|u| \ge T^{1/2}}\mcg(\sigma+it+iu)K_{\eta}(u)|R(t)|^{2}\Phi\Big(\frac{t\log T}{T}\Big)\mathrm{d}u\mathrm{d}t, \\
    W_3(T) &\, \coloneqq\int_{|t| \ge T \log T}\int_\mbr\mcg(\sigma + it + iu)K_{\eta}(u)|R(t)|^{2}\Phi\Big(\frac{t\log T}{T}\Big)\mathrm{d}u\mathrm{d}t.
\end{align*}
The following upper bound for $W_2(T)$ can be established by using Lemma \ref{approximatededekind} and the fact that $K_{\eta}(u) \ll u^{-4\phi(d)}$
\begin{align}
\label{W2Tupperbound}
    W_2(T) &\, \ll d^{\phi(d)} \int_{1\le |t| \le T \log T}\int_{|u| \ge T^{1/2}} (1 + |t| + |u|)^{3\phi(d)/10} K_{\eta}(u)|R(t)|^{2}\Phi\Big(\frac{t\log T}{T}\Big)\mathrm{d}u\mathrm{d}t  \nonumber \\    
    &\, \ll d^{\phi(d)}(T\log T)^{3\phi(d)/10} T^{(3\phi(d)/10 - 4\phi(d)+1)/2} \int_\mbr|R(t)|^{2}\Phi\Big(\frac{t\log T}{T}\Big)\mathrm{d}t \nonumber \\
    &\, = o\Big(\frac{T}{\log T} \sum_{n \in \mbn} f(n)^2\Big).
\end{align}
Furthermore, for $W_3(T)$, applying Lemma \ref{approximatededekind}, we obtain
$$
W_3(T) \ll d^{\phi(d)} \int_{|t| \ge T\log T}|t|^{3\phi(d)/10}|R(t)|^{2}\Phi\Big(\frac{t\log T}{T}\Big)\mathrm{d}t \int_\mbr |u|^{3\phi(d)/10} K_{\eta}(u)\mathrm{d}u.
$$
Since $\Phi$ decays rapidly, 
\begin{align}
\label{W3Tupperbound}
    W_3(T) &\, \ll d^{\phi(d)} (T\log T)^{3\phi(d)/10} \Phi \Big(\frac{(\log T)^2}{\sqrt{2}}\Big) \int_\mbr |R(t)|^{2} \Phi \Big(\frac{t\log T}{\sqrt{2}T}\Big)\mathrm{d}t \nonumber \\
    & \, = o \Big(\frac{T}{\log T} \sum_{n \in \mbn} f(n)^2\Big).
\end{align}
\eqref{W2Tupperbound} and \eqref{W3Tupperbound} yield the following result:
\begin{equation}
    \label{WTW1To}
    W(T) = W_1(T) + o\Big(\frac{T}{\log T} \sum_{n \in \mbn} f(n)^2\Big) \ll (\widehat{K_{\eta}}(0)\mcz(\sigma)+o(1))\frac{T}{\log T}\sum_{n \in \mbn} f(n)^2,
\end{equation}
where
$$
W(T) \coloneqq \int_{|t|\ge1}\int_\mbr\mcg(\sigma + it + iu)K_{\eta}(u)|R(t)|^{2}\Phi\Big(\frac{t\log T}{T}\Big)\mathrm{d}u\mathrm{d}t.
$$
\par
Then, define
\begin{align*}
    &W_0(T) \coloneqq \int_{|t| \ge 1}\Big(\sum_{n=1}^{\infty}\frac{\widehat{K_{\eta}}(\log n)a(n)}{n^{\sigma +it}}\Big)|R(t)|^{2}\Phi\Big(\frac{t\log T}{T}\Big)\mathrm{d}t, \\
    &\mce \coloneqq \int_{|t| \ge1} 2\pi \tau \cdot |R(t)|^{2}\Phi\Big(\frac{t\log T}{T}\Big)\mathrm{d}t.
\end{align*}
Cauchy's integral formula implies $\mce = o(T^{\lambda} \sum_{n \in \mbn}f(n)^2)$, thus,
\begin{equation}
    \label{WT=W0T-GAMMA}
    W(T) = W_0(T) + o\Big(T^{\lambda} \sum_{n \in \mbn}f(n)^2\Big).
\end{equation}
\par
To extend the $t$ integral to $\mbr$, we consider the integral over the region $|t| \le 1$. Plainly, 
$$
\mce^{\prime}\coloneqq\int_{|t| \le 1}\Big(\sum_{n=1}^{\infty}\frac{\widehat{K_{\eta}}(\log n)a(n)}{n^{\sigma +it}}\Big)|R(t)|^{2}\Phi\Big(\frac{t\log T}{T}\Big)\mathrm{d}t \ll R(0)^{2} \sum_{n=1}^{\infty}\frac{\widehat{K_{\eta}}(\log n)a(n)}{n^{\sigma}}.
$$
Combining with Lemma \ref{Ketauproperties}, we deduce that
\begin{equation}
\label{theregiontxiaoyu1}
    \mce^{\prime}\ll R(0)^{2} \widehat{K_{\eta}}(0) (\log_2 T)^A \sum_{n\leqslant T^{2\varepsilon}}1 \ll \widehat{K_{\eta}}(0) T^{\lambda + 2 \varepsilon} \sum_{n \in \mbn}f(n)^2.
\end{equation}
Recalling that $\lambda + 2\varepsilon <1$ and combining \eqref{WTW1To}, \eqref{WT=W0T-GAMMA} and \eqref{theregiontxiaoyu1}, it follows that
\begin{equation}
\label{ZKsigmaTgg}
    \mcz(\sigma) \gg \int_\mbr\Big(\sum_{n=1}^{\infty}\frac{\widehat{K_{\eta}}(\log n)a_{\mathbb{K},\ell}(n)}{n^{\sigma +it}}\Big)|R(t)|^{2}\Phi\Big(\frac{t\log T}{T}\Big)\mathrm{d}t + o(1).
\end{equation}
\par
Denoting $\widetilde{W}(T)$ as the integral on the right-side of the formula above, we have
$$
\widetilde{W}(T)= \sqrt{2 \pi} \frac{T}{\log T} \sum_{m,n \in \mcm^{\prime}_d} \sum_{k=1}^{\infty}\frac{\widehat{K_{\eta}}(\log n)a_\mbk(k)r(m)r(n)}{k^{\sigma}} \Phi\Big(\frac{T}{\log T} \log \frac{km}{n}\Big).
$$
Here we note the fact that $a(n) \ge a_\mbk(n)$. Following the idea of \cite{bondarenko2023dichotomy}, we restrict the inner sum to $k \le T^{\varepsilon/3\eta}$, yielding
$$
\widetilde{W}(T) \gg \widehat{K_{\eta}}(0) \frac{T}{\log T} \sum_{m,n\in\mathcal{M}^{\prime}_d} \sum_{k\le T^{\varepsilon/3n}} \frac{a(k)r(m)r(n)}{k^{\sigma}} \Phi\Big(\frac{T}{\log T} \log \frac{km}{n}\Big).
$$
Employing the techniques similar to those in \cite[Eq. (21)]{bondarenko2018argument}, we obtain
\begin{align}
\label{widetildeWTlowerbound}
\widetilde{W}(T)
&\, \gg \widehat{K_{\eta}}(0) \frac{T}{\log T}\sum_{m,n\in\mcm_d} \sum_{\substack{km = n \\k\le T^{\varepsilon/3\eta}}} \frac{a_\mbk(k)f(m)f(n)}{k^{\sigma}} \nonumber \\
&\, = \widehat{K_{\eta}}(0) \frac{T}{\log T} \sum_{n\in\mcm_d} \frac{f(n)}{n^{\sigma}} \sum_{q\mid n, q\ge n/T^{\varepsilon/3\eta}} a_\mbk\Big(\frac{n}{q}\Big)f(q)q^{\sigma}.
\end{align}
Combining Propositions \ref{prop41}, \ref{prop42}, and \ref{prop43} yields
$$\mcz(\sigma) \gg \exp\bigg((\delta\gamma c_d + o(1))\frac{(\log T  )^{1-\sigma}(\log_3 T)^{\sigma}}{(\log_2 T)^{\sigma}}\bigg) + o(1)$$
for $\lambda + 2\varepsilon <1.$ Then, 
\begin{align*}
    \max_{-T \log T - T^{1/2}\le t\le T \log T +T^{1/2}}&|\zeta_\mbk^{(\ell)}(\sigma+it)| \gg \mcz(\sigma) + O(1)\\
    &\gg \exp\bigg((\delta\gamma c_d + o(1))\frac{(\log T)^{1-\sigma} (\log_3 T)^{\sigma}}{(\log_2 T)^{\sigma}}\bigg) + o(1).
\end{align*}
Let $\lambda +2 \varepsilon <1$ with small $\varepsilon$. By choosing $\delta$, $\gamma$ and $c_d$ to be sufficiently close, respectively, to $1$, $1$ and $\sqrt{\phi(d)/(e-1)}$, and applying the same transform as in Section \ref{proofthm1.1process}, we conclude the proof of Theorem \ref{thm1.2}.

\section{Supplementary Conclusions}
\label{supplementaryconclusions}
In Sections \ref{proofofthm1.1} and \ref{proofofthm1.2}, we present Theorems \ref{thm1.1} and \ref{thm1.2} along with their proofs. These theorems establish lower bounds for maximum of derivatives of $\zeta_\mbk(s)$ on and near the critical line.
\par
Inspired by \cite{dong2023Onde}, we present the following conclusion. Since the method used is almost identical to that in Section \ref{proofofthm1.2}, we omit the proof here and recommend \cite{dong2025largevalue} for further details.
\begin{thm}
    \label{thm5.1}
    Let $A$ and $D$ be arbitrary positive numbers, and $T$ be sufficiently large. Set $\sigma_D \coloneqq 1/2 + D/(\log_2 T)$. Then uniformly for $d \ll (\log_2T)^A,$ we have
    $$
    \max_{0\le t\le T}|\zeta^{(\ell)}_\mbk(\sigma_D+it)| \ge \exp\bigg(\Big(\frac{\sqrt{\phi(d)}}{\sqrt{(e-1)}\cdot e^D} + o(1)\Big)\sqrt{\frac{\log T\log_3 T}{\log_2 T}}\bigg).
    $$   
\end{thm}
Unfortunately, the set $\mcm_d$ and the resonator $R(t)$ constructed in Section \ref{resonator2} can only be used to study extreme values within the region where $|\sigma-1/2| \ll (\log_2T)\inv$. For the larger range beyond this, we supplement the following result.
\begin{thm}
     \label{thm5.2}
     Let $A$ be an arbitrary positive number. If $T$ is sufficiently large in the range $1/2 \le \sigma \le 2/3$, then uniformly for $d \ll (\log_2T)^A,$ we have
$$
    \max_{t\in[0,T]}|\zeta_\mbk^{(\ell)}(\sigma+it)|\ge \exp\bigg(\Big(\sqrt{\frac{\phi(d)}{e-1}}+o(1)\Big)\frac{(\log T)^{1-\sigma}}{(\log_2 T)^{\sigma}}\bigg).
$$
\end{thm}
\begin{proof}
In fact, we only need to redefine the set $\mcm_{k,d}$ as
$$\mcm_{k,d} \coloneqq \bigg\{ n \in \mathrm{supp}(f) : n \text{ has at least } \Delta_k \coloneqq \frac{\alpha (\log N)^{2 - 2\sigma}}{k^2 \log_3 N} \text{ prime divisors in } \mcp_{k,d} \bigg\},
$$
and the value of $f(p)$ for $p \in \mcp_d$ as 
$$
f(p) = c_d \bigg( \frac{(\log N)^{1 - \sigma} (\log_2 N)^{\sigma}}{\log_3 N} \bigg) \frac{1}{p^{\sigma} (\log p - \log_2 N - \log_3 N - \log \phi(d))}.
$$
Here, the definitions of $\mcp_d$ and $\mcp_{k,d}$ remain unchanged, and the remainder of the argument proceeds as in Section \ref{proofthm1.2process}.
\end{proof}
\par
It is worth noting that using a similar method in \cite{bondarenko2018note}, we can extend the above result almost to the region $1/2 < \sigma < \sigma_0$ for any fixed $\sigma_0$ such that $1/2<\sigma_0<1$
\begin{equation}
    \label{finalbound}
        \max_{t\in[0,T]}|\zeta_\mbk^{(\ell)}(\sigma+it)|\ge \exp\bigg(C(\sigma) \phi(d)^{1-\sigma}\frac{(\log T)^{1-\sigma}}{(\log_2 T)^{\sigma}}\bigg).
\end{equation}
The constant $C(\sigma)$ depends on $\sigma$. Thus, in the critical strip, derivatives of $\zeta_\mbk(s)$ attain extreme values of the same order of magnitude as those in \cite[Theorem 2.(B)]{yang2022extreme} for $\zeta(s)$.
\par
It can be seen from \eqref{finalbound} that the closer $\sigma$ is to the 1-line, the worse this lower bound becomes. When $\sigma = 1$, to the best of our knowledge, an effective method to obtain a sharp lower bound has not yet been established. For $\zeta(s)$, Levinson \cite{levinson2017Omega} showed that when $\sigma=1$
\begin{align*}
 \max_{t\in[1,T]} |\zeta(1+it)| \ge e^E \log_2 T +O(1), 
\end{align*}
where $E$ is the Euler-Mascheroni constant\footnote{The best known bound for extreme values of $\zeta(s)$ on the 1-line can be found in \cite{aistleitner2019extreme}.}. For $\sigma=1$, better results on derivatives of $\zeta(s)$ can be found in \cite{yang2024extreme}. For $\sigma$ near the 1-line, see also \cite{dong2023Onde}. 

\appendix
\section{Proof of Proposition \ref{agcdlowerbound}}
In this appendix, we provide a supplementary proof of Proposition \ref{agcdlowerbound}. It is easy to show that 
\begin{equation}
    \label{smakprod}
    S_{1/2}(\mcm,a_\mbk) = \prod_{1 \le k \le (\log_2 N)^\delta}S_{1/2}(\mcm_k,a_\mbk),
\end{equation}
thus we consider $S_{1/2}(\mcm_k,a_\mbk)$ for a fixed $k$ first. Following the argument of \cite[Lemma 2]{fonga2023extreme}, we obtain 
\begin{equation}
\label{61lowerbound}
    S_{1/2}(\mcm_k,a_\mbk) \ge \sum_{\substack{\ell,\ell^\prime \mid W\\ \omega(\ell),\omega(\ell^\prime)\le W_k/2}} \frac{(\ell,\ell^\prime)a_\mbk^\prime(\ell)a_\mbk^\prime(\ell^\prime)}{a_\mbk^\prime((\ell,\ell^\prime))^2\sqrt{\ell\ell^\prime}} \sum_{\substack{q,q^\prime \mid W \\ (q,\ell)=(q^\prime,\ell^\prime)=1 \\\omega(q),\omega(q^\prime)\le W_k/2}}\frac{(q,q^\prime)a_\mbk(q)a_\mbk(q^\prime)}{a_\mbk^\prime((q,q^\prime))^2\sqrt{qq^\prime}},
\end{equation}
where the multiplicative function $a_\mbk^\prime(n)$ is defined by
$$
a_\mbk^\prime(n) = \prod_{p \mid n}\frac{\phi(d)+1}{2}.
$$
It is clear that $a_\mbk^\prime(p) = (\phi(d)+1)/2 \ge \phi(d)/2$ for all primes $p$.
\par
Let $\mcf_k \coloneqq \mcf_k(\ell,\ell^\prime)$ denote the inner sum in \eqref{61lowerbound}. Set $q=n_1d_1$ and $q^\prime=n_1^\prime d_1$, respectively. Then choosing $n=(q,q^\prime)$ in the identity $n=\sum_{d_1 \mid n}\phi(d_1)$, it follows that 
$$
\mcf_k =\sum_{\substack{d_1 \mid W \\ \omega(d_1) \leq W_k / 2 \\(d_1, \ell \ell^{\prime})=1}} \frac{\phi(d_1)}{d_1} \sum_{\substack{n_1, n_1^{\prime} \mid W \\(n_1, \ell d_1)=(n_1^{\prime}, \ell^{\prime} d_1)=1 \\ \omega(n_1), \omega(n_1^{\prime}) \leqslant w_k / 2-\omega(d_1)}} \frac{a_\mbk(n_1) a_\mbk(n_1^{\prime})}{a_\mbk((n_1, n_1^{\prime}))^2 \sqrt{n_1 n_1^{\prime}}}.
$$
Since every term in $\mcf_k$ is positive, it suffices to consider the case $(n_1,n_1^\prime)=1$, which yields that
$$
\mcf_k \ge  \sum_{\substack{d_1 \mid W \\
\omega(d_1) \leq W_k / 2 \\
(d_1, \ell \ell^{\prime})=1}} \frac{\phi(d_1)}{d_1}\sum_{\substack{n_1^{\prime} \mid W \\
(n_1^{\prime}, \ell^{\prime} d_1)=1 \\
\omega(n_1^{\prime}) \leq W_k / 2-\omega(d_1)}} \frac{a_\mbk(n_1^{\prime})}{\sqrt{n_1^{\prime}}} \sigma\Big(a_\mbk, \frac{W_k}{2}-\omega(d_1), \ell d_1 n_1^{\prime}\Big),
$$
where 
$$
\sigma(a_\mbk, R, r)\coloneqq\sum_{\substack{n \mid W \\
\omega(n) \leq R \\
(n, r)=1}} \frac{a_\mbk(n)}{\sqrt{n}}.
$$
Note that $\sigma(a_\mbk, R, r)$ increases as $R$ increases. Since $\phi(d_1)/d_1 \gg 1$ as $d_1\mid W$, we have
$$
\mcf_k \gg \sum_{\substack{d_1 \mid W \\ \omega(d_1) \leq W_k / 2 \\(d_1, \ell \ell^{\prime})=1}} \sum_{\substack{n_1^{\prime} \mid W \\(n_1^{\prime}, \ell^{\prime} d_1)=1 \\ \omega(n_1^{\prime}) \leq W_k / 2-\omega(d_1)}} \frac{a(n_1^{\prime})}{\sqrt{n_1^{\prime}}} \sigma\Big(a_\mbk, \frac{W_k}{2} -\omega(d_1), \ell d_1 n_1^{\prime}\Big).
$$
\par
Set $\ell=n_2d_2$ and $\ell^\prime=n_2^\prime d_2$, respectively. Applying the same calculation to the outer sum in \eqref{61lowerbound}, we obtian
\begin{align*}
\label{62lowerbound}
   S_{1/2}(\mcm_k,a_\mbk) \gg &\,
   \sum_{\substack{d_1,d_2\mid W \\\omega(d_1)\le W_k/2  \\ \omega(d_2) \le W_k/2 \\(d_1,d_2)=1 }} \sum_{\substack{n_2,n_2^\prime \mid W\\(n_2n_2^\prime,d_1d_2)=1\\\omega(n_2)\le W_k/2-\omega(d_2)\\\omega(n_2^\prime)\le W_k/2-\omega(d_2)\\(n_2,n_2^\prime)=1}}\frac{a_\mbk^\prime(n_2)a_\mbk^\prime(n_2^\prime)}{\sqrt{n_2n_2^\prime}} \\
   &\, \times \sum_{\substack{n_1^{\prime} \mid W \\(n_1^{\prime}, n_2^\prime d_2d_1)=1 \\ \omega(n_1^{\prime}) \leq W_k / 2-\omega(d_1)}} \frac{a(n_1^{\prime})}{\sqrt{n_1^{\prime}}} \sigma\Big(a_\mbk, \frac{W_k}{2} -\omega(d_1),n_2 d_2 d_1 n_1^{\prime}\Big).
\end{align*}
\par
Next, we estimate the lower bound for the above. Set
$$
u_k \coloneqq \bigg\lfloor\frac{\nu}{k}\sqrt{\frac{\log N}{\log_2 N \log_3 N}} \bigg\rfloor,
$$
where $\nu$ is a bounded parameter. Restrict the previous outer sum to pairs $(d_1,d_2)$ such that $\omega(d_1),\omega(d_2)\le W_k/2-u_k$, and the inner sum to integers $n_2,n_2^\prime$ such that $\omega(n_2)=\omega(n_2^\prime)=u_k$. Combining with the fact that
$$
\sigma\Big(a_\mbk, \frac{W_k}{2}-\omega\left(d_1\right), n_2 d_2 d_1 n_1^{\prime}\Big) \geq \sigma(a_\mbk, u_k, n_2 d_2 d_1 n_1^{\prime}),
$$
we have 
\begin{align*}
    & S_{1/2}(\mcm_k, a_\mbk) \gg \\
    & \sum_{\substack{d_1, d_2 \mid W \\ \omega(d_1) \le W_k / 2-u_k \\ \omega(d_2) \le W_k / 2-u_k\\(d_1, d_2)=1}} \sum_{\substack{n_2, n_2^{\prime} \mid W \\(n_2 n_2^{\prime}, d_1 d_2)=1 \\ \omega(n_2), \omega(n_2^{\prime})=u_k \\(n_2, n_2^{\prime})=1}} \frac{a_\mbk^{\prime}\left(n_2\right) a_\mbk^{\prime}\left(n_2^{\prime}\right)}{\sqrt{n_2^{\prime} n_2}} \sum_{\substack{n_1^{\prime} \mid W \\(n_1^{\prime}, d_2 n_2^{\prime} d_1)=1 \\ \omega(n_1^{\prime}) \le W_k / 2-\omega(d_1)}}\frac{a_\mbk(n_1^{\prime})}{\sqrt{n_1^{\prime}}} \sigma(a_\mbk, u_k, n_2 d_2 d_1 n_1^{\prime}).
\end{align*}
Proceeding as in the calculations of \cite[pp. 8-9]{tenen2019galsum}, we obtain 
\begin{equation}
    \label{Hupperbound}
    \sigma(a_\mbk, u_k, n_2 d_2 d_1 n_1^{\prime}) \ge H^{u_k}e^{o(u_k)},
\end{equation}
where 
$$
H \coloneqq \frac{2 k e\sqrt{\phi(d)}(\sqrt{\alpha}-1) \alpha^{k / 2}}{\nu} \sqrt{\log _3 N}\bigg(1+O\bigg(\sqrt{\frac{\log _2 N \log _3 N}{\log N}}\bigg)\bigg).
$$
In fact, the Siegel-Walfisz theorem implies 
\begin{equation}
    \label{sumpk}
    \sum_{p \in P_k} \frac{a_\mbk(p)}{\sqrt{p}} = 2\sqrt{\phi(d)}\alpha^{k/2}(\sqrt{\alpha}-1)\sqrt{\frac{\log N}{\log_2 N}}\Big(1+O\Big(\frac{k+\log\phi(d)+\log_3 N}{\log_2 N} \Big) \Big).
\end{equation}
Furthermore, noting that $\omega(n_1^\prime d_1)+\omega(n_2 d_2) \le W_k$, we have 
\begin{equation}
    \label{sumpkmid}
    \sum_{\substack{p \in P_k \\ p \mid n_2 n_1^\prime d_2 d_1}}\frac{a_\mbk(p)}{\sqrt{p}} \le \frac{\sqrt{\phi(d)}W_k}{\alpha^{k/2}\sqrt{\log N \log_2 N}} 
\end{equation} 
Thus \eqref{sumpk} and \eqref{sumpkmid} yields that
\begin{equation}
    \sum_{\substack{p \in P_k \\ p \mid n_2 n_1^\prime d_2 d_1}}\frac{a_\mbk(p)}{\sqrt{p}}  = 2\sqrt{\phi(d)}\alpha^{k/2}(\sqrt{\alpha}-1)\sqrt{\frac{\log N}{\log_2 N}}\Big(1+O\Big(\frac{1}{\log_3 N} \Big) \Big).
\end{equation}
Together with the fact that 
$$
\sigma(a_\mbk, u_k, n_2 d_2 d_1 n_1^{\prime}) \gg \frac{1}{w_k!}\bigg(\sum_{\substack{p \in P_k \\ p \nmid n_2 n_1^\prime d_2 d_1}}\frac{a_\mbk(p)}{\sqrt{p}}\bigg)^{w_k}
$$
and the Stirling's formula, we obtain \eqref{Hupperbound}. A similar calculation yields 
$$
\sum_{\substack{n_1^{\prime} \mid N_k \\\left(n_1^{\prime}, d_2 n_2^{\prime} d_1\right)=1 \\ \omega(n_1^{\prime}) \leq W_k / 2-\omega(d_1)}} \frac{a_\mbk(n_1^{\prime})}{\sqrt{n_1^{\prime}}}=\sigma\Big(a_\mbk, \frac{W_k}{2}-\omega(d_1), n_2^{\prime} d_2 d_1\Big) \gg H^{w_k} e^{o(w_k)},
$$
therefore we have the following lower bound
\begin{equation}
    \label{smkaklower1}
    S_{1/2}(\mcm_k,a_\mbk) \gg H^{2w_k}e^{o(w_k)} \sum_{\substack{d_1,d_2 \mid W \\ \omega(d_1),\omega(d_2) \le W_k/2 \\ (d_1,d_2)=1}} \sum_{\substack{n_2,n_2^\prime \mid W \\ (n_2n_2^\prime,d_1d_2)=1 \\ \omega(n_2),\omega(n_2^\prime)\le W_k/2-\omega(d_2)\\(n_2,n_2^\prime)=1}} \frac{a_\mbk^\prime(n_2)a_\mbk^\prime(n_2^\prime)}{\sqrt{n_2n_2^\prime}}.
\end{equation}
By performing a similar calculation again and noting that $a_\mbk^\prime(p) = (\phi(d)+1)/2 \ge \phi(d)/2$, we obtain 
$$
\sum_{\substack{n_2,n_2^\prime \mid W \\ (n_2n_2^\prime,d_1d_2)=1 \\ \omega(n_2),\omega(n_2^\prime)\le W_k/2-\omega(d_2)\\(n_2,n_2^\prime)=1}} \frac{a_\mbk^\prime(n_2)a_\mbk^\prime(n_2^\prime)}{\sqrt{n_2n_2^\prime}} \gg \Big(\frac{H}{2} \Big)^{2w_k}e^{o(w_k)}.
$$
Substituting the above into \eqref{smkaklower1} shows that
\begin{equation}
    \label{smkaklower2}
    S_{1/2}(\mcm_k,a_\mbk) \gg \Big(\frac{H}{\sqrt{2}} \Big)^{4w_k}e^{o(w_k)} \sum_{\substack{d_1,d_2 \mid W \\ \omega(d_1),\omega(d_2) \le W_k/2-w_k \\ (d_1,d_2)=1}} 1 \eqqcolon \Big(\frac{H}{\sqrt{2}} \Big)^{4w_k}e^{o(w_k)}V_k.
\end{equation}
\par
Using the same calculation as \cite[Eq. (2.13)]{tenen2019galsum}, we have
$$
V_k \ge |\mcm_k|\Big(\frac{\beta+o(1)}{2k^2 \alpha^k(\alpha-1)\log_3 N} \Big)^{2w_k}.
$$
Thus combining \eqref{smkaklower2} we find
$$
\frac{S_{1/2}(\mcm_k,a_\mbk)}{|\mcm_k|} \gg \Big(\frac{h}{v^2} \Big)^{2w_k}e^{o(w_k)},
$$
where $h \coloneqq  e^2 \beta \phi(d)(\sqrt{\alpha}-1)/(\sqrt{\alpha}+1)$.
\par
Using \eqref{smakprod} and the definition of $w_k$, we have
$$
\frac{S_{1/2}(\mcm, a_\mbk)}{|\mcm|} \ge \exp \bigg(\big(\rho+o(1)\big)\sqrt{\frac{\log N \log_3 N}{\log_2 N}}\bigg),
$$
where $\rho \coloneqq 2 \delta \nu \log(h/\nu^2)$.
Setting $v = \sqrt{h}/e$ gives 
$$
\rho = \frac{4\delta\sqrt{h}}{e} = 4 \delta \sqrt{\frac{\beta\phi(d)(\sqrt{\alpha}-1)}{\sqrt{\alpha}+1}}.
$$
By choosing $\beta\delta\log\alpha \to 1$ and $\alpha \to 1$, we finish the proof of Proposition \ref{agcdlowerbound}.

\section*{Acknowledgments}
The authors would like to thank Dr. Winston Heap for his corrections to the manuscript. The third author is supported by the Natural Science Foundation of Henan Province (Grant No. 252300421782).

	\bibliographystyle{siam}
    \bibliography{reference}
\end{document}